\documentclass[12pt,fleqn,leqno]{article}

\author{
V.\ Oproiu \\
\and D.D.\ Poro\c sniuc\\  }
\date{}
\title{A K\"ahler Einstein structure on the cotangent
bundle of a Riemannian manifold \thanks{partially supported by the
Grant 100/2003, Ministerul Educa\c tiei \c si Cercet\u arii,
Rom\^ania}}

\begin{document}

\maketitle
\begin{abstract}
We use the natural lifts of the fundamental tensor field $g$ to
the cotangent bundle $T^*M$ of a Riemannian manifold $(M,g)$, in
order to construct an almost Hermitian structure $(G,J)$ of
diagonal type on $T^*M$. The obtained almost complex structure $J$
on $T^*M$ is integrable if and only if the base manifold has
constant sectional curvature and the second coefficient, involved
in its definition is expressed as a rational function of the first
coefficient and its first order derivative. Next one shows that
the obtained almost Hermitian structure is almost K\"ahlerian.
Combining the obtained results we get a family of K\"ahlerian
structures on $T^*M$, depending on one essential parameter. Next
we study the conditions under which the considered K\"ahlerian
structure is Einstein. In this case $(T^*M,G,J)$ has constant
holomorphic curvature.


\end{abstract}

\vskip5mm {\large \bf Introduction} \vskip5mm

 The differential geometry of the cotangent bundle $T^*M$ of a
 Riemannian manifold $(M,g)$ is quite similar to that of the
 tangent bundle $TM$. However there are some differences, due to
 the fact that the lifts (vertical, complete, horizontal etc.) to
$T^*M$ cannot be defined just like in the case of $TM$.

In the present paper we study a family of natural K\"ahler
Einstein structures $(G,J)$, of diagonal type induced on $T^*M$
from the Riemannian metric $g$. They are obtained in a manner
quite similar to that used in \cite{Oproiu5} (see also
\cite{Oproiu4}) by using a similar parametrization. The considered
natural Riemannian metric $G$ of diagonal type on $T^*M$ is
defined by using two parameters $u,v$ which are smooth functions
depending on the energy density $t$ on $T^*M$. The vertical and
horizontal distributions are orthogonal to each other and the dot
products induced on them from $G$ are isomorphic (isometric) by
duality.

Next, the family of the natural almost complex structures $J$ on
$T^*M$ that interchange the vertical and horizontal distributions
depends on the same two essential parameters $u, v$. From the
integrability condition for $J$ it follows that the base manifold
$M$ must have constant sectional curvature $c$ and the second
parameter $v$ must be expressed as a rational function depending
on the first parameter $u$ and its derivative. Of course, in the
obtained formula there are involved too the constant $c$ and the
energy density $t$.

 Next it follows that $G$ is Hermitian with respect to $J$ and it
 follows that the fundamental $2$-form $\phi$, associated to the
 almost Hermitian structure $(G,J)$ is the fundamental form defining the
 usual symplectic structure on $T^*M$, hence it is closed. In the case
 where the integrability condition for $J$ is
fulfilled, we get a K\"ahlerian structure on $T^*M$ and this
structure depends on one essential parameter $u$.

In the case where the considered K\"ahlerian structure is Einstein
we get a second order differential equation fulfilled by the
parameter $u$ and we have been able to find the general solution
of this equation. We have a generic case, when the curvature $c$
is negative and whole $(T^*M,G,J)$ is K\"ahler Einstein. Then
there is another case where $c$ is positive and $(G,J)$ is defined
on a tube around the zero section in $T^*M$. In both cases one
obtains that $(T^*M,G,J)$ has constant holomorphic curvature. The
cases where some singularities occur will be discussed in some
forthcoming papers.

The manifolds, tensor fields and geometric objects we consider in
this paper, are assumed to be differentiable of class $C^{\infty}$
(i.e. smooth). We use the computations in local coordinates but
many results from this paper may be expressed in an invariant
form. The well known summation convention is used throughout this
paper, the range for the indices $h,i,j,k,l,r,s$ being always
${\{}1,...,n{\}}$ (see \cite{GheOpr}, \cite{OprPap1},
\cite{OprPap2}, \cite{YanoIsh}). We shall denote by
${\Gamma}(T^*M)$ the module of smooth vector fields on $T^*M$.

\vskip5mm {\large \bf 1. Some geometric properties of $T^*M$}
\vskip5mm

Let $(M,g)$ be a smooth $n$-dimensional Riemannian manifold and
denote its cotangent bundle by $\pi :T^*M\longrightarrow M$.
Recall that there is a structure of a $2n$-dimensional smooth
manifold on $T^*M$, induced from the structure of smooth
$n$-dimensional manifold  of $M$. From every local chart
$(U,\varphi )=(U,x^1,\dots ,x^n)$  on $M$, it is induced a local
chart $(\pi^{-1}(U),\Phi )=(\pi^{-1}(U),q^1,\dots , q^n,$
$p_1,\dots ,p_n)$ on $T^*M$ as follows. For a cotangent vector
$p\in \pi^{-1}(U)\subset T^*M$, the first $n$ local coordinates
$q^1,\dots ,q^n$ are  the local coordinates $x^1,\dots ,x^n$ of
its base point $x=\pi (p)$ in the local chart $(U,\varphi )$ (in
fact we have $q^i=\pi^* x^i=x^i\circ \pi, \ i=1,\dots n)$. The
last $n$ local coordinates $p_1,\dots ,p_n$ of $p\in \pi^{-1}(U)$
are the vector space coordinates of $p$ with respect to the
natural basis $(dx^1_{\pi(p)},\dots , dx^n_{\pi(p)})$, defined by
the local chart $(U,\varphi )$,\ i.e. $p=p_idx^i_{\pi(p)}$. Due to
this special structure of differentiable manifold for $T^*M$ it is
possible to introduce the concept of $M$-tensor field on it. An
$M$-tensor field of type $(r,s)$ on $T^*M$ is defined by sets of
$n^{r+s}$ components (functions depending on $q^i$ and $p_i$),
with $r$ upper indices and $s$ lower indices, assigned to induced
local charts $(\pi^{-1}(U),\Phi )$ on $T^*M$, such that the local
coordinate change rule is that of the local coordinate components
of a tensor field of type $(r,s)$ on the base manifold $M$ (see
\cite{Mok} for further details in the case of the tangent bundle);
e.g., the components $p_i,\ i=1,\dots ,n$, corresponding to the
last $n$ local coordinates of a cotangent vector $p$, assigned to
an induced local chart $(\pi^{-1}(U), \Phi )$ define an $M$-tensor
field of type $(0,1)$ on $T^*M$. A usual tensor field of type
$(r,s)$ on $M$ may be thought of as an $M$-tensor field of type
$(r,s)$ on $T^*M$. If the considered tensor field on $M$ is
covariant only, the corresponding $M$-tensor field on $T^*M$ may
be identified with the induced (pullback by $\pi $) tensor field
on $T^*M$. Some useful $M$-tensor fields on $T^*M$ may be obtained
as follows. Let $u:[0,\infty ) \longrightarrow {\bf R}$ be a
smooth function and let $\|p\|^2=g^{-1}_{\pi(p)}(p,p)$ be the
square of the norm of the cotangent vector $p\in \pi^{-1}(U)$
($g^{-1}$ is the tensor field of type (2,0) on $M$ having as
components the entries $g^{ij}(x)$ of the inverse of the matrix
$(g_{ij}(x))$ defined by the components of $g$ in the local chart
$(U,\varphi )$). If $\delta ^i_j$ are the Kronecker symbols (in
fact, they are the local coordinate components of the identity
tensor field $I$ on $M$), then the components $u(\|p\|^2)\delta
^i_j$ define an $M$-tensor field of type $(1,1)$ on $T^*M$.
Similarly, if $g_{ij}(x)$ are the local coordinate components of
the metric tensor field $g$ on $M$ in the local chart $(U,\varphi
)$, then the components $u(\|p\|^2) g_{ij}(\pi(p))$ define a
symmetric
 $M$-tensor field of type $(0,2)$ on $T^*M$. The components
 $g^{0i}=p_hg^{hi}$, as well as $u(\|p\|^2)g^{0i}$, define $M$-tensor
 fields of type $(1,0)$ on $T^*M$. Of course,  all the components
 considered above are in the induced local chart $(\pi^{-1}(U),\Phi)$.

We shall use the horizontal distribution $HT^*M$, defined by the
Levi Civita connection $\dot \nabla $ of $g$, in order to define
some first order natural lifts to $T^*M$ of the Riemannian metric
$g$ on $M$. Denote by $VT^*M= {\rm Ker}\ \pi _*\subset TT^*M$ the
vertical distribution on $T^*M$. Then we have the direct sum
decomposition

\begin{equation}
TT^*M=VT^*M\oplus HT^*M.
\end{equation}

If $(\pi^{-1}(U),\Phi)=(\pi^{-1}(U),q^1,\dots ,q^n,p_1,\dots
,p_n)$ is a local chart on $T^*M$, induced from the local chart
$(U,\varphi )= (U,x^1,\dots ,x^n)$, the local vector fields
$\frac{\partial}{\partial p_1}, \dots , \frac{\partial}{\partial
p_n}$ on $\pi^{-1}(U)$ define a local frame for $VT^*M$ over $\pi
^{-1}(U)$ and the local vector fields $\frac{\delta}{\delta
q^1},\dots ,\frac{\delta}{\delta q^n}$ define a local frame for
$HT^*M$ over $\pi^{-1}(U)$, where
$$
\frac{\delta}{\delta q^i}=\frac{\partial}{\partial
q^i}+\Gamma^0_{ih} \frac{\partial}{\partial p_h},\ \ \ \Gamma
^0_{ih}=p_k\Gamma ^k_{ih}
 $$
and $\Gamma ^k_{ih}(\pi(p))$ are the Christoffel symbols of $g$.

The set of vector fields $(\frac{\partial}{\partial p_1},\dots
,\frac{\partial}{\partial p_n}, \frac{\delta}{\delta q^1},\dots
,\frac{\delta}{\delta q^n})$ defines a local frame on $T^*M$,
adapted to the direct sum decomposition (1). Remark that
$$
\frac{\partial}{\partial p_i}=(dx^i)^V,\ \ \frac{\delta}{\delta
q^i}=(\frac{\partial}{\partial x^i})^H,
$$
where $\theta^V$ denotes the vertical lift to $T^*M$ of the
$1$-form $\theta$ on $M$ and $X^H$ denotes the horizontal lift to
$T^*M$ of the vector field $X$ on $M$.

Now we shall present the following auxiliary result. \vskip5mm

{\bf Lemma 1}. \it If $n>1$ and $u,v$ are smooth functions on
$T^*M$ such that
$$
u g_{ij}+v p_ip_j=0,\ p\in \pi^{-1}(U),
$$
on the domain of any induced local chart on $T^*M$, then $u=0,\
v=0$.\rm \vskip 0.5cm

The proof is obtained easily by transvecting the given relation
with components $g^{ij}$ of the tensor field $g^{-1}$ and $g^{0j}$
(Recall that the functions $g^{ij}(x)$ are the components of the
inverse of the matrix $(g_{ij}(x))$, associated to $g$ in the
local chart $(U,\varphi )$ on $M$).\vskip5mm

\bf Remark. \rm From the relations of the type
$$
u g^{ij}+v g^{0i}g^{0j}=0,\ p\in \pi^{-1}(U),
$$
$$
u\delta ^i_j+vg^{0i} p_j=0,\ p\in \pi^{-1}(U),
$$
it is obtained, in a similar way, $u=v=0$. We have used the
notation $g^{0i}=p_hg^{hi}$. \vskip5mm

Since we work in a fixed local chart $(U,\varphi )$ on $M$ and in
the corresponding induced local chart $(\pi^{-1}(U),\Phi )$ on
$T^*M$, we shall use the following simpler (but less clear)
notations
$$
\frac{\partial}{\partial p_i}=\partial ^i,\ \ \frac{\delta}{\delta
q^i}= \delta _i.
$$

Denote by
\begin{equation}
t=\frac{1}{2}\|p\|^2=\frac{1}{2}g^{-1}_{\pi(p)}(p,p)=\frac{1}{2}g^{ik}(x)p_ip_k,
\ \ \ p\in \pi^{-1}(U)
\end{equation}
the energy density defined by $g$ in the cotangent vector $p$. We
have $t\in [0,\infty)$ for all $p\in T^*M$. For a vector field $X$
on $M$ we shall denote by $g_X$ the $1$-form on $M$ defined by
$g_X(Y)=g(X,Y)$, for all vector fields $Y$ on $M$. For a $1$-form
$\theta$ on $M$, we shall denote by $\theta^\sharp=g^{-1}_\theta$
the vector field on $M$ defined by the usual musical isomorphism,
i.e. $g(\theta^\sharp,Y)=\theta (Y)$, for all vector field $Y$ on
$M$. Remark that, for $p\in T^*M$, we can consider the vector
$p^\sharp$, tangent to $M$ in $\pi(p)$.

\vskip5mm {\large \bf 2. A natural almost K\"ahler structure of
diagonal type on the cotangent bundle} \vskip5mm

From now on we shall work in a fixed local chart $(U,\varphi)$ on
$M$ and in the induced local chart $(\pi^{-1}(U),\Phi)$ on $T^*M$.
Consider two real valued smooth functions $u,v$ defined on
$[0,\infty)\subset {\bf R}$ and define the following $M$-tensor
field of type $(0,2)$ defined by the Riemannian metric $g$ on the
cotangent bundle $T^*M$ of the $n$-dimensional Riemannian manifold
$(M,g)$
\begin{equation}
G_{ij}(p)= u(t) g_{ij}(\pi(p))+v(t)p_ip_j.
\end{equation}

It follows easily that the matrix $(G_{ij})$ is positive definite
if and only if $u>0,\ u+2tv>0$. The inverse of this matrix has the
entries

\begin{equation}
H^{kl}(p)= \frac{1}{u(t)} g^{kl}(\pi(p))+w(t)g^{0k}g^{0l},
\end{equation}
where $g^{0k}=p_hg^{hk}$ and
\begin{equation}
w= -\ \frac{v}{u(u+2tv)}.
\end{equation}

One shows easily that the components $H^{kl}$ assigned to the
induced local chart $(\pi^{-1}(U),q^1,\dots, q^n,p_1,\dots ,p_n)$
on $T^*M$ define an $M$-tensor field of type $(2,0)$. The property
of the matrix $(H^{kl})$ to be the inverse of the matrix
$(G_{ij})$ is expressed by the formulas
$$
G_{ik}H^{kj}=\delta_i^j.
$$

{\bf Remark.} If the matrix $(G_{ij})$ is positive definite, its
inverse $(H^{kl})$ is positive definite too. This aspect can be
seen directly since we have
$$
\frac{1}{u}>0,\ \ \frac{1}{u}+2tw =\frac{1}{u+2tv}>0.
$$

Using the $M$-tensor fields defined by $G_{ij},\ H^{kl}$, the
following Riemannian metric may be considered on $T^*M$
\begin{equation}
G=G_{ij}dq^idq^j+H^{ij}Dp_iDp_j,
\end{equation}
where $Dp_i=dp_i-\Gamma^0_{ij}dq^j$ is the absolute (covariant)
differential of $p_i$ with respect to the Levi Civita connection
$\dot\nabla$ of $g$ (We have used the notation
$\Gamma^0_{ij}=p_h\Gamma^h_{ij}$, where $\Gamma^h_{ij}$ are the
Christoffel symbols defined by $g$). Equivalently, we have
$$
G(\delta_i,\delta_j)=G_{ij},~~~G(\partial^i ,\partial^j)=H^{ij},~~
G(\partial^i,\delta_j)=~ G(\delta_j,\partial^i)=0.
$$

Remark that $HT^*M,~VT^*M$ are orthogonal to each other with
respect to $G$, but the Riemannian metrics induced from $G$ on
$HT^*M,~VT^*M$ are not the same, so the considered metric $G$ on
$T^*M$ is not a metric of Sasaki type. However, the matrix
associated to the dot product induced from $G$ on $VT^*M$ is the
inverse of the matrix associated to the dot product induced from
$G$ on $HT^*M$, so that the metrics induced from $G$ on $VT^*M, \
HT^*M$ could be considered as being isomorphic (isometric). The
$2n\times 2n$-matrix associated to $G$, with respect to the
adapted local frame $(\frac{\delta}{\delta q^1},\dots
,\frac{\delta}{\delta q^n},\frac{\partial}{\partial p_1},\dots
,\frac{\partial}{\partial p_n})$ has two $n\times n$-blocks on the
first diagonal
\begin{displaymath}
G= \left(
\begin{array}{cc}
G_{ij} & 0  \\
0 & H^{ij}
\end{array}
\right).
\end{displaymath}

The Riemannian metric $G$ is called a \it natural lift of diagonal
type \rm of $g$. Remark also that the system of 1-forms
$(dq^1,...,dq^n,Dp_1,...,Dp_n)$ defines a local frame on
$T^{*}T^*M$, dual to the local frame $(\delta_1 ,...,\delta_n,~
\partial^1,...,\partial^n)$ adapted to the
direct sum decomposition (1).

Next, an almost complex structure $J$ is defined on $T^*M$ by the
same $M$-tensor fields $G_{ij},\ H^{kl}$, expressed in adapted
local frames by
$$
J{\delta _i}=G_{ik}\partial ^k,\ \ J{\partial ^i}=-H^{ik}\delta
_k.
$$

The matrix of $J$ with respect to the adapted local basis
$(\delta_1 ,...,\delta_n,~ \partial^1,...,\partial^n)$ is
\begin{displaymath}
J= \left(
\begin{array}{cc}
0 &  -H^{ij}  \\
G_{ij} & 0
\end{array}
\right).
\end{displaymath}

From the property of the $M$-tensor field $H^{kl}$ to be defined
by the inverse of the matrix defined by the components of the
$M$-tensor field $G_{ij}$, it follows easily that $J$ defines an
almost complex structure on $T^*M$. The almost complex structure
defined by $J$ is called a \it natural lift of diagonal type \rm
of $g$. \vskip2mm

{\bf Proposition 2}. \it The total space of the cotangent bundle
$T^*M$, endowed with the Riemannian metric $G$ and the almost
complex structure $J$ (both natural lifts of $g$ of diagonal type)
has a structure of almost K\"ahlerian manifold. \vskip2mm

 Proof. \rm Since the matrix $(H^{kl})$ is the inverse of the
 matrix $(G_{ij})$, it follows easily that
 $$
G(J\delta_i,J\delta_j)=G(\delta_i,\delta_j)=G_{ij}, \
G(J\partial^i, J\partial^j)=G(\partial^i,\partial^j)=H^{ij},
$$
$$
G(J\delta_i,J\partial^j)=G(J\partial^j,J\delta_i)=G(\delta_i,\partial^j)=
G(\partial^j,\delta_i)=0
$$

Hence $G(JX,JY)=G(X,Y)$, for all vector fields $X,Y$ on $T^*M$.
Thus $(T^*M,G,J)$ is an almost Hermitian manifold. The fundamental
$2$-form associated with this almost Hermitian structure is
$\phi$, defined by
$$
\phi(X,Y) = G(X,JY),
$$
for all vector fields $X,Y$ on $T^*M$. By a straightforward
computation we get
$$
\phi(\delta_i,\delta_j)=G(\delta_i,G_{jk}\partial^k)=0, \
\phi(\partial^i,\partial^j)=G(\partial^i,H^{jk}\delta_k)=0,
$$
$$
-\phi(\delta_j,\partial^i)=\phi(\partial^i,\delta_j)=G(\partial^i,
G_{jk}\partial^k)=G_{jk}H^{ik}= \delta^i_j
$$

It follows that
$$
\phi =Dp_i\wedge dq^i= dp_i\wedge dq^i,
$$
due to the symmetry of $\Gamma^0_{ij}=p_h\Gamma^h_{ij}$. It
follows that $\phi$ does coincide with the fundamental $2$-form
defining the usual symplectic structure on $T^*M$. Of course, we
have $d\phi =0$, i.e. $\phi$ is closed. Therefore $(T^*M,G,J)$ is
an almost K\"ahler manifold. \pagebreak

\vskip5mm {\large \bf 3. A natural K\"ahler structure on $T^*M$}
\vskip5mm

We shall study the integrability of the almost complex structure
defined by $J$ on $T^*M$. To do this we need the following well
known formulas for the brackets of the vector fields
$\partial^i=\frac{\partial}{\partial p_i},\delta_i=
\frac{\delta}{\delta q^i},~ i=1,...,n$
\begin{equation}
[\partial^i,\partial^j]=0;~~~[\partial^i,\delta_j]=\Gamma^i_{jk}\partial^k;~~~
[\delta_i,\delta_j] =R^0_{kij}\partial^k,
\end{equation}
where $R^h_{kij}(\pi(p))$ are the local coordinate components of
the curvature tensor field of $\dot \nabla$ on $M$ and
$R^0_{kij}(p)=p_hR^h_{kij}$ . Of course, the components
 $R^0_{kij}$, $R^h_{kij}$ define M-tensor fields of types
 (0,3), (1,3) on $T^*M$, respectively. \vskip2mm

{\bf Theorem 3. } {\it The Nijenhuis tensor field of the almost
complex structure $J$ on $T^*M$ is given by}

$$
\left\{
\begin{array}{l}
N(\delta_i,\delta_j)=\{(v(2tu^\prime-u)+uu^\prime)(\delta^h_ig_{jk}-
\delta^h_jg_{ik})-R^h_{kij}\}p_h\partial^k,
\\ \mbox{ } \\
N(\delta_i,\partial^j)=H^{kl}H^{jr}\{(v(2tu^\prime-u)+uu^\prime)(\delta^h_ig_{rl}-
\delta^h_rg_{il})-R^h_{lir}\}p_h\delta_k,
\\ \mbox{ } \\
N(\partial^i,\partial^j)=H^{ir}H^{jl}\{(v(2tu^\prime-u)+uu^\prime)(\delta^h_lg_{rk}-
\delta^h_rg_{lk})-R^h_{klr}\}p_h\partial^k.
\end{array}
\right.
$$
\vskip2mm

{\it Proof. } Recall that the Nijenhuis tensor field $N$ defined
by $J$ is given by
$$
N(X,Y)=[JX,JY]-J[JX,Y]-J[X,JY]-[X,Y],\ \ \forall X,Y \in \Gamma
(T^*M).
$$
Then, we have $\delta_kt =0,\ \partial^k t = g^{0k}$ and
$\dot\nabla_iG_{jk}=0,\ \dot\nabla_iH^{jk}= 0$, where
$$
\dot\nabla_iG_{jk}=
\delta_iG_{jk}-\Gamma^l_{ij}G_{lk}-\Gamma^l_{ik}G_{lj}
$$
$$
\dot\nabla_iH^{jk}=
\delta_iH^{jk}+\Gamma^j_{il}H^{lk}+\Gamma^k_{il}H^{lj}
$$

The above expressions for the components of $N$ can be obtained by
a quite long, straightforward  computation.\vskip2mm

{\bf Theorem 4} {\it The almost complex structure $J$ on $T^*M$ is
integrable if and only if the base manifold $M$ has constant
sectional curvature $c$ and the function $v$ is given by}

\begin{equation}
v=\frac{c-uu^\prime}{2tu^\prime-u}.
\end{equation}

{\it Proof.} From the condition $N=0$ one obtains
$$
\{(v(2tu^\prime-u)+uu^\prime)(\delta^h_ig_{jk}-
\delta^h_jg_{ik})-R^h_{kij}\}p_h=0.
$$
Differentiating with respect to $p_l$, taking $p_h=0~\forall h \in
{\{}1,...,n{\}}$, it follows that the curvature tensor field of
$\dot \nabla$ has the expression

$$
R^l_{kij}=u(0)[-v(0)+u^\prime(0)]({\delta}^l_ig_{jk}-{\delta}^l_jg_{ik}).
$$
Using the Schur theorem (in the case where $M$ is connected and
$dim~M \geq 3$), it follows that $(M,g)$ has the constant
sectional curvature $c=u(0)[-v(0)+u^\prime(0)]$. Then we obtain
the expression (8) of $v$.

 Conversely, if $(M,g)$ has
constant sectional curvature $c$ and $v$ is given by (8), it
follows in a straightforward way that $N=0$.\vskip2mm

\bf Remark. \rm In the case where $u^2-2ct=0$, we have $uu^\prime-
c=0,\ u-2tu^\prime=0$ too. So, this case must be thought of as a
singular case and should be considered separately. Recall that the
function $u$ must fulfill the conditions
$$
u>0,\ u+2tv=\frac{2ct-u^2}{2tu^\prime-u}>0.
$$

Thus the family of natural K\"ahlerian structures of diagonal type
on $T^*M$ (when $N=0$) depends on one essential coefficient $u$
satisfying some the supplementary conditions.

\vskip5mm {\large \bf 4. The Levi Civita connection of the metric
$G$ and its curvature tensor field} \vskip5mm

The Levi Civita connection $\dot\nabla$ on a Riemannian manifold
$(M,g)$ is determined by the conditions
$$
\dot\nabla g=0,~~~~~\dot T=0,
$$
where $\dot T$ is its torsion tensor field. The explicit
expression of this connection is obtained from the formula
$$
2g({\dot\nabla}_XY,Z)=X(g(Y,Z))+Y(g(X,Z))-Z(g(X,Y))+
$$
$$
+g([X,Y],Z)-g([X,Z],Y)-g([Y,Z],X), ~~~~~\forall
X,Y,Z{\in}{\Gamma}(M).
$$

We shall use this formula in order to obtain the expression of the
Levi Civita connection ${\nabla}$ of $G$ on $T^*M$. The final
result can be stated as follows \vskip2mm
\newpage
{\bf Theorem 7.} {\it The Levi Civita connection ${\nabla}$ of $G$
on $T^*M$ has the following expression in the local adapted frame
$(\partial^1,...,\partial^n,~ \delta_1,..., \delta_n)$

$$
\left\{
\begin{array}{l}
\nabla_{\partial^i}\partial^j =Q^{ij}_h\partial h,\ \ \ \ \ \
\nabla_{\delta_i}\partial^j=-\Gamma^j_{ih}\partial^h+P^{hj}_i\delta_h,
\\ \mbox{ } \\
\nabla_{\partial^i}\delta_j=P^{hi}_j\delta_h,\ \ \ \ \ \ \
\nabla_{\delta_i}\delta_j=\Gamma^h_{ij}\delta_h+S_{hij}\partial^h,
\end{array}
\right.
$$
where $Q^{ij}_h, P^{hi}_j, S_{hij}$ are $M$-tensor fields on
$T^*M$, defined by
$$
\left\{
\begin{array}{l}
Q^{ij}_h = \frac{1}{2}G_{hk}(\partial^iH^{jk}+
\partial^jH^{ik} -\partial^kH^{ij}),
\\ \mbox{ } \\
P^{hi}_j=\frac{1}{2}H^{hk}(\partial^iG_{jk}-H^{il}R^0_{ljk}),
\\ \mbox{ } \\
S_{hij}=-\frac{1}{2}G_{hk}\partial^kG_{ij}+\frac{1}{2}R^0_{hij}.
\end{array}
\right.
$$ \rm

After replacing of the expressions of the involved $M$-tensor
fields and their derivatives, one obtains

$$
\left\{
\begin{array}{l}
Q^{ij}_h =
-\frac{u^\prime}{2u}(\delta^i_hg^{0j}+\delta^j_hg^{0i})-
\frac{v(u^\prime+2u^2w)}{2u^3w}g^{ij}p_h- \frac{v(2u^\prime
w+uw^\prime)}{2u^2w}g^{0i}g^{0j}p_h ,
\\ \mbox{ } \\
P^{hi}_j=\frac{u^\prime}{2u}\delta^h_jg^{0i}-\frac{(c+uv)w}{2v}
\delta^i_jg^{0h}+\frac{uv-c}{2u^2}g^{ih}p_j+\frac{vw(uv-c)+uw(u^\prime
v-uv^\prime)}{2uv}g^{0i}g^{0h}p_j,
\\ \mbox{ } \\
S_{hij}= \frac{c-uv}{2}g_{jh}p_i-\frac{c+uv}{2}g_{ih}p_j+
\frac{u^\prime v}{2uw}g_{ij}p_h+\frac{v(v^\prime
-2uvw)}{2uw}p_hp_ip_j.
\end{array}
\right.
$$

In the case of a K\"ahler structure on $T^*M$, the final
expressions of these $M$-tensor fields can be obtained by doing
the necessary replacements of the functions $v,w$ from (8) and (5)
and their derivatives. However, the final expressions are quite
complicate but they may be obtained quite automatically by using
the Mathematica package RICCI for doing tensor computations (see
\cite{Lee}}).

Now we shall indicate the obtaining of the components of the
curvature tensor field of the connection $\nabla$.

The curvature tensor $K$ field of the connection $\nabla $ is
obtained from the well known formula
$$
K(X,Y)Z=\nabla_X\nabla_YZ-\nabla_Y\nabla_XZ-\nabla_{[X,Y]}Z.
$$

The components of $K$ with respect to the adapted local frame
$(\partial^1,\dots ,\partial^n,$ $\delta_1,\dots ,\delta_n)$ are
obtained easily:
\newpage

$$
\left\{
\begin{array}{l}
K(\partial^i,\partial^j)\partial^k=PPP^{ijk}_h\partial^h
=(\partial^iQ^{jk}_h-\partial^jQ^{ik}_h+Q^{jk}_lQ^{il}_h-
Q^{ik}_lQ^{jl}_h)\partial^h,
\\ \mbox{ } \\
K(\partial^i,\partial^j)\delta_k=PPQ^{ijh}_k\delta_h
=(\partial^iP^{hj}_k-\partial^jP^{hi}_k+P^{lj}_kP^{hi}_l-
P^{li}_kP^{hj}_l)\delta_h,
\\ \mbox{ } \\
K(\delta_i,\delta_j)\partial^k=QQP^k_{ijh}\partial^h = (-R^k_{hij}
- R^0_{lij}Q^{lk}_h+S_{hil}P^{lk}_j-S_{hjl}P^{lk}_i)\partial^h,
\\ \mbox{ } \\
K(\delta_i,\delta_j)\delta_k=QQQ_{ijk}^h\delta_h
=(R^h_{kij}-R^0_{lij}P^{hl}_k+
S_{ljk}P^{hl}_i-S_{lik}P^{hl}_j)\delta_h,
\\ \mbox{ } \\
K(\partial^i,\delta_j)\delta_k=PQQ^i_{jkh}\partial^h =(\partial^i
S_{hjk}+S_{ljk}Q^{il}_h-S_{hjl}P^{li}_k)\partial^h,
\\ \mbox{ } \\
K(\partial^i, \delta_j)\partial^k=PQP_j^{ikh}\delta_h =(\partial^i
P^{hk}_j+P^{hi}_lP^{lk}_j-Q^{ik}_lP^{hl}_j)\delta_h.
\end{array}
\right.
$$

 The explicit expressions of these components are obtained after
some quite long and hard computations, made by using the package
RICCI.

Next, the components of the Ricci tensor field are obtained as
traces of $K$:
$$
\left\{
\begin{array}{l}
Ric(\partial^j,\partial^k)=RicPP^{jk}=PPP^{hjk}_h-PQP^{jkh}_h,
\\ \mbox{ } \\
Ric(\delta_j,\delta_k)=RicQQ_{jk}=QQQ_{hjk}^h+PQQ^h_{jkh},
\\ \mbox{ } \\
Ric(\partial^j,\delta_k)=Ric(\delta_k,\partial^j)=0.
\end{array}
\right.
$$

\vskip5mm {\large \bf 5 The cotangent bundle $T^*M$ as a K\"ahler
\break Einstein manifold} \vskip5mm

Doing the necessary computations, we obtain the final expressions
of the components of the Ricci tensor field of $\nabla$
$$
\left\{
\begin{array}{l}
RicQQ_{jk}=\frac{a}{2(u-2tu^\prime)^2}g_{jk}+\frac{\alpha}{2u^2(u-2tu^\prime)^4}p_jp_k,
\\ \mbox{ } \\
RicPP^{jk}=\frac{a}{2u^2(u-2tu^\prime)^2}g^{jk}+\frac{\beta}{2u^2(u^2-2ct)(u-2tu^\prime)^2}
g^{0j}g^{0k},
\end{array}
\right.
$$
where the coefficients $a,\alpha,\beta $ are given by
$$
a=n(u-2tu^\prime)(2cu-2ctu^\prime-u^2u^\prime)+2(2ct-u^2)(tuu^{\prime\prime}+
uu^\prime-tu^{\prime 2}),
$$
$$
\alpha =n(u-2tu^\prime)(-2c^2u^3+6c^2tu^2u^\prime+3cu^4u^\prime-12
c^2t^2uu^{\prime 2}-3 u^5u^{\prime 2}+8c^2 t^3u^{\prime 3} -
$$
$$
-4c t^2u^2u^{\prime 3}+ 4tu^4u^{\prime 3}
-4c^2t^2u^2u^{\prime\prime}+
4ctu^4u^{\prime\prime}-u^6u^{\prime\prime})+
$$\pagebreak
$$
+2(2ct-u^2)(-3cu^3u^\prime+7ctu^2u^{\prime 2}+4u^4u^{\prime 2} -8c
t^2uu^{\prime 3}-8tu^3u^{\prime 3}+4ct^3u^{\prime
4}+4t^2u^2u^{\prime 4}-
$$
$$
-7ctu^3u^{\prime\prime}+2u^5u^{\prime\prime} +6ct^2u^2u^\prime
u^{\prime\prime}+3tu^4u^\prime u^{\prime\prime}- 6t^2u^3u^{\prime
2}u^{\prime\prime}-8ct^3u^2u^{\prime\prime 2}+
$$
$$
+4t^2u^4u^{\prime\prime 2}- 2ct^2u^3u^{(3)}
+tu^5u^{(3)}+4ct^3u^2u^\prime u^{(3)}- 2t^2u^4u^\prime u^{(3)}),
$$

$$
\beta =n(u-2tu^\prime)(2c^2u-2c^2tu^\prime-3cu^2u^\prime+
6ctuu^{\prime 2}-u^3u^{\prime 2}-4ct^2u^{\prime 3}+2tu^2u^{\prime
3}+
$$
$$
+2ctu^2u^{\prime\prime}-u^4u^{\prime\prime})+
$$
$$
+2(2c^2tuu^\prime +cu^3u^\prime-2c^2t^2u^{\prime
2}-5ctu^2u^{\prime 2}-2u^4u^{\prime 2}+8ct^2uu^{\prime
3}+4tu^3u^{\prime 3}- 4ct^3u^{\prime 4}-
$$
$$
-2t^2u^2u^{\prime 4}+2c^2t^2uu^{\prime\prime}+
5ctu^3u^{\prime\prime}-2u^5u^{\prime\prime}-6ct^2u^2u^\prime
u^{\prime\prime}-tu^4u^\prime u^{\prime\prime}+ 4t^2u^3u^{\prime
2}u^{\prime\prime}+
$$
$$
+8ct^3u^2u^{\prime\prime 2}-4t^2u^4 u^{\prime\prime
2}+2ct^2u^3u^{(3)}-tu^5u^{(3)}-4ct^3u^2u^\prime
u^{(3)}+2t^2u^4u^\prime u^{(3)}).
$$

 In order to find the conditions under which $(T^*M,G,J)$ is Einstein,
 we consider the differences
$$
\left\{
\begin{array}{l}
DiffQQ_{jk}=RicQQ_{jk}-\frac{a}{2u(u-2tu^\prime)^2}\ G_{jk},
\\ \mbox{ } \\
DiffPP^{jk}=RicPP^{jk}-\frac{a}{2u(u-2tu^\prime)^2}\ H^{jk},
\end{array}
\right.
$$

whose explicit expressions are
$$
\left\{
\begin{array}{l}
DiffQQ_{jk}=\frac{u^2-2ct}{2u^2(u-2tu^\prime)^4}\ \gamma \ p_jp_k,
\\ \mbox{ } \\
DiffPP^{jk}=\frac{1}{2u^2(u^2-2ct)(u-2tu^\prime)^2}\ \gamma \
g^{0j}g^{0k}.
\end{array}
\right.
$$
The expression of the factor $\gamma$ is given by
$$
\gamma =n(u^2-2ct)(2tu^\prime-u)(u^2u^{\prime\prime} -2tu^{\prime
3}+2uu^{\prime 2})+
$$
$$
+2(2cu^3u^\prime-4ctu^2u^{\prime 2}-3u^4u^{\prime 2}+
6ct^2uu^{\prime 3}+5tu^3u^{\prime 3}-4ct^3u^{\prime 4}-
2t^2u^2u^{\prime 4}+
$$
$$
+6ctu^3u^{\prime\prime}-2u^5u^{\prime\prime}-4ct^2u^2u^\prime
u^{\prime\prime}-2tu^4u^\prime u^{\prime\prime}+ 4t^2u^3u^{\prime
2}u^{\prime\prime}+8ct^3u^2u^{\prime\prime 2}-
$$
$$
-4t^2u^4u^{\prime\prime 2}+2ct^2u^3u^{(3)}-tu^5u^{(3)}-
4ct^3u^2u^\prime u^{(3)}+2t^2u^4u^\prime u^{(3)}).
$$

We are interested in finding the functions $u$ for which
$DiffQQ_{jk}=0$, $DiffPP^{jk}=0$. We should exclude the case
$u^2-2ct=0$ which leads to a singularity. Thus we must see what
happens in the case $\gamma=0, \ u^2-2ct\neq 0$. Generally, this
equation is almost impossible to solve. However, it is reasonable
to ask for the solution $u$ to be independent of the dimension $n$
of $M$. In this case we must have

\begin{equation}
u^2u^{\prime\prime} -2tu^{\prime 3}+2uu^{\prime 2}=0
\end{equation}

The general solution of this equation is obtained transforming it
in an equation in the inverse function $t=t(u)$ of $u$. One
obtains an equation of Euler type from which we get

$$
u=A\pm\sqrt{A^2+Bt},\ \ \ A>0,\ B\neq 0,\ A^2+Bt>0,
$$
where $A,B$ are the integration constants.

\bf Remark. \rm The equation $(9)$ has two other singular
solutions
$$
u=A,\ \ \ A>0,\ A^2-2ct>0,
$$
$$
u=At,\ \ \ A>0,\ 2c-A^2t>0,
$$
which will be discussed in some forthcoming papers.

From now on we shall consider only the case of the function
$u=A+\sqrt{A^2+Bt}$. Asking for the found function to be a
solution for the remaining part of the equation $\gamma =0$, one
finds $B=-2c$. Thus the general solution for the condition for
$(T^*M,G,J)$ to be Einstein is

\begin{equation}
u = A+\sqrt{A^2-2ct}, \ \ \ A>0,\ A^2-2ct>0.
\end{equation}

If the constant sectional curvature $c$ of $M$ is negative, the
solution $u$ is defined on the whole $T^*M$.

 If $c$ is positive,
the solution $u$ is defined only in the tube around the zero
section in $T^*M$, defined by $0\leq \|p\|^2<\frac{A^2}{c}$. Then
we obtain easily
$$
v= \frac{1}{2t}(A-\frac{4ct}{A}-\sqrt{A^2-2ct}),\ \ \
w=\frac{-A^3+3Act+(A^2-2ct)^{3/2}}{4ct^2(A^2-2ct)}.
$$
Remark that $v, w$ are well defined even if $t=0$. Then we have
$$
u+2tv=\frac{2(A^2-2ct)}{A}>0,
$$
thus the conditions for the existence of the K\"ahler Einstein
manifold $(T^*M,G,J)$ are fulfilled.

The components of the Ricci tensor field defined by $G$ are
$$
RicQQ_{jk}=\frac{(n+1)c}{A}G_{jk},\
RicPP^{jk}=\frac{(n+1)c}{A}H^{jk},
$$
$$
Ric(\partial^j,\delta_k)=Ric(\delta_k,\partial^j)=0.
$$

Hence we may state our main result. \vskip2mm

{\bf Theorem 5.} {\it 1. Assume that the Riemannian manifold
$(M,g)$ has constant negative sectional curvature $c$. Then
$(T^*M,G,J)$ defined by $u$ given in (10), with $A>0$ is a
K\"ahler Einstein manifold.

2. Assume that $(M,g)$ has constant positive sectional curvature
$c$. Then the tube around the zero section in $T^*M$, defined by
the condition $0\leq\|p\|^2<\frac{A^2}{c}$ is a K\"ahler Einstein
manifold, with $(G,J)$ defined as in the case 1.}

\vskip5mm {\large \bf 6 The holomorphic sectional curvature of
\break $(T^*M,G,J)$} \vskip5mm

Recall that a K\"ahlerian manifold $(M,g,J)$ has constant
holomorphic sectional curvature $k$ if its curvature tensor field
$R$ is given by
$$
R(X,Y)Z=\frac{k}{4}(g(Z,Y)X-g(Z,X)Y+
$$
$$
+g(Z,JY)JX-g(Z,JX)JY+ 2g(X,JY)JZ),
$$
for all vector fields $X,Y,Z$ defined on $M$.

In the case of the K\"ahler Einstein structure $(G,J)$ on $T^*M$
(on a tube around zero section in $T^*M$) obtained in Theorem 5,
we can check by a straightforward computation that the components
of the curvature tensor field $K$ of $\nabla$ are given by
$$
K(\delta_i,\delta_j)\delta_k=\frac{c}{2A}(\delta^h_iG_{jk}-
\delta^h_jG_{ik})\delta_h,\ \
K(\delta_i,\delta_j)\partial^k=\frac{c}{2A}(\delta^k_jG_{ih}-
\delta^k_iG_{jh})\partial^h,
$$
$$
K(\partial^i,\partial^j)\delta_k=\frac{c}{2A}(\delta^j_kH^{ih}-
\delta^i_kH^{jh})\delta_h,\ \
K(\partial^i,\partial^j)\partial^k=\frac{c}{2A}(\delta^i_hH^{jk}-
\delta^j_hH^{ik})\partial^h,
$$
$$
K(\partial^i,\delta_j)\delta_k=\frac{c}{2A}(\delta^i_hG_{jk}+
\delta^i_kG_{jh}+2\delta^i_jG_{kh})\partial^h,
$$
$$
K(\partial^i,\delta_j)\partial^k=-\frac{c}{2A}(\delta^h_jH^{ik}+
\delta^k_jH^{ih}+2\delta^i_jH^{kh})\delta_ h,
$$

From these relations we get \vskip2mm

{\bf Theorem 6.} {\it The K\"ahler Einstein structure $(G,J)$ on
$T^*M$ (on a tube around zero section in $T^*M$) obtained in
Theorem 5 has constant holomorphic sectional curvature
$k=\frac{2c}{A}$.}

\vskip 1.5cm
\begin{minipage}{2.5in}
\begin{flushleft}
V.Oproiu\\
Faculty of Mathematics\\
University "Al.I.Cuza", Ia\c si \\
Rom\^ania.\\
e-mail: voproiu@uaic.ro
\end{flushleft}
\end{minipage}
\hfill
\begin{minipage}{2.5in}
\begin{flushleft}
D.D.Poro\c sniuc\\
Department of Mathematics\\
National College "M. Eminescu" \\
Boto\c sani, Rom\^ania.\\
e-mail: danielporosniuc@lme.ro
\end{flushleft}
\end{minipage}
\end{document}